\font\title=cmbx10 at 12pt
\magnification=1200

\normalbaselineskip = 18pt
\noindent

\normalbaselines
\parskip=10pt plus1pt minus.5pt

\centerline{\title  A UNIFORM KADEC-KLEE PROPERTY FOR}
\centerline {\title SYMMETRIC OPERATOR SPACES}
 \smallskip
\centerline {{\sl P.G. Dodds, T.K. Dodds*,
 P.N. Dowling, C.J. Lennard*$\dag$ and F.A. Sukochev}
 \footnote{$\ $}
 {* Research partially supported by A.R.C. \hfil\break
  \indent $\dag$ Research partially supported by a U. Pittsburgh FAS Grant
 \hfil\break
 }}

\medskip

\font\title=cmbx10 at 12pt

\magnification = 1200
\normalbaselineskip = 20pt
\parskip = 12pt plus 1pt minus 0.5pt
\normalbaselines

\def\p{\,\prec\kern-1.2ex \prec\,}
\def\ch{\raise 0.5ex \hbox{$\chi$}}

\def\nm{{\cal M}}

\def\nmt{\widetilde {\cal M}}

\def\r{{\bf R}^+}
\def\jr{{\bf R}}
\def\jc{{\bf C}}
\def\jn{{\bf N}}

\font\title=cmbx10 at 12pt

\noindent {\bf ABSTRACT}

\smallskip

We show that if a rearrangement invariant Banach function
space $E$ on the positive semi-axis satisfies a
non-trivial lower $q-$ estimate with constant $1$ then the
corresponding space $E(\nm)$ of $\tau-$measurable operators,
affiliated with
 an arbitrary semi-finite von Neumann algebra $\nm$ equipped with a
distinguished faithful, normal, semi-finite trace $\tau $,
 has the uniform Kadec-Klee
 property for the topology of local convergence in measure.  In particular, the
Lorentz function spaces $L_{q,p}$
 and the Lorentz-Schatten classes ${\cal C}_{q,p}$ have the UKK property
 for convergence locally in measure and for the weak-operator topology,
 respectively. As a partial converse , we show that if $E$ has the UKK property
with respect
 to local convergence in measure then $E$ must satisfy some non-trivial
 lower $q$-estimate. We also prove a uniform Kadec-Klee result for
 local convergence in any
 Banach lattice satisfying a lower $q$-estimate.

\smallskip

\noindent {\bf 1. Introduction and preliminaries}

\smallskip
It is well-known (see, for example, [Ar],[Si]) that the Schatten $p-$classes
${\cal C}_p, 1\leq p<\infty ,$ have the (so-called ) Kadec-Klee (or Radon-Riesz)
property, that is norm convergence in ${\cal C}_p$  coincides with weak operator
convergence for sequences in the unit sphere. It has been shown recently by
Lennard [Le1] that the direct argument given by Arazy [Ar] for trace-class
operators may be refined to show that the Schatten class ${\cal C}_1$ has
the uniform Kadec-Klee property for the weak operator
topology (see below for the precise definitions). Further extensions to the
setting of trace ideals have been announced by Sukochev [Su].  From a somewhat
different (though related) viewpoint, it is a classical theorem of F.Riesz that
norm convergence for sequences in the
unit sphere of $L^1 ([0,1])$ coincides with
convergence in measure. Appropriate uniform versions of this theorem may be
found in [Le2], and extensions to certain Lorentz spaces in [Le3].

It is the intention of the
present paper to give simultaneous extensions of each of the above theorems in
the setting of symmetric spaces of measurable operators. Our approach is based
on a  notion of local convergence in measure for measurable operators which
specialises to the weak operator topology for the case of trace ideals in the
sense of Gohberg and Krein [GK], to pointwise convergence in the case of
symmetric sequence spaces, and to convergence in measure on sets of finite
measure in the case of rearrangement invariant function spaces. Let $\nm $ be a
semi-finite von Neumann algebra, equipped with a distinguished faithful, normal
semi-finite trace $\tau $. Our principal result (Theorem 2.7) shows that if $E$
is a symmetric (rearrangement-invariant) quasi-Banach function space on the
semi-axis $\r$ which satisfies a non-trivial lower estimate with constant $1$,
then the corresponding symmetric operator space $E(\nm )$ has the uniform
Kadec-Klee property for local convergence in measure. Moreover, our theorem is
best possible in the sense that symmetric spaces on the positive half-line which
have this uniform property must necessarily satisfy a non-trivial lower
estimate.   The existence of such an estimate is a geometric
condition which is readily verified, and so our results essentially give a
useful geometric characterisation of this uniform  Kadec-Klee type property
which not only unifies the results alluded to above, but also yields an
extensive class of new examples. In each of these examples, our principal
theorem yields results concerning normal structure and fixed point
properties. For a more detailed discussion of these aspects and the relation
between normal structure and uniform Kadec-Klee type properties, we refer the
interested reader to [Le2], and the references contained therein. 

We remark that Dilworth and Hsu [DH] have recently proved a uniform version
 of the theorem of [Ar]/[Si] i.e. they show that if a symmetric
 sequence space $E$ has the UKK property for pointwise convergence,
 then the associated unitary matrix space ${\cal C}_E$ has the UKK
 property for the weak operator topology. They complement this Banach space
 result by showing that the Schatten class ${\cal C}_p$ has 
 UKK(weak operator) for all $ 0 < p \le 1$, which is a special case of
 our Theorem 2.7 below. Along the way they prove a special case of 
 our main inequality in Proposition 2.5. However, their methods are based on a 
 generalization of McCarthy's characterization of Schatten class norms
 in [M],
 rather than inequalities for decreasing rearrangements of
 operators; and so differ
 from the methods we have used.

We begin by recalling some necessary terminology.  Let $E$ be a complex quasi-
Banach lattice.
If $0<p<\infty $, then  $E$ is said to be $p-${\it convex} (resp.
$p-${\it concave }) if there exists a constant $C>0$ such
that for all finite sequences $(x_n)$ in $E$\hfil \break
\centerline {$\Vert (\sum |x_n|^p)^{1\over p}\Vert
 _{_{E}}\leq C(\sum \Vert x_n\Vert^p_{_{E}})^{1\over p},$}
\centerline {({\rm resp.}$\  C\Vert(\sum
|x_n|^p)^{1\over p}\Vert _{_{E}}\geq (\sum\Vert x_n\Vert
^p_{_{E}})^{1\over p }$).}
\smallskip
The least such constant $C$ is called the $p-${\it
convexity } (resp. $p-${\it concavity}) constant of $E$ and
is denoted by $M^{(p)}\left (E\right )$ (resp.
$M_{(p)}\left (E\right)).$  For $0<p<\infty ,\  E^{(p)}$
will denote the quasi-Banach lattice defined by
$$ E^{(p)}=\{x:|x|^p\in E\}$$
equipped with the quasi-norm
$$\Vert x\Vert _{_{E^{(p)}}}=\Vert |x|^p\Vert
_{_{E}}^{1/p}.$$

It is not difficult to see that if $E$ is $\alpha-$convex
and $q-$concave then $E^{(p)}$ is $\alpha p-$convex and
$qp-$concave with $M^{(\alpha p)}\left (E^{(p)}\right )\leq
M^{(\alpha )}\left ( E\right )^{1\over p }$ and $ M_{(qp)}\left
(E^{(p)}\right )\leq M_{(q)}\left (E \right )^{1\over p}.$
Consequently, if $E$ is $\alpha-$convex then $E^{(1/\alpha
)}$ is $1-$convex and so can be renormed as a Banach
lattice (cf.[LT]).

The quasi-Banach lattice $E$ is said to satisfy a
$lower\ q-estimate$ (resp. $upper\ p-estimate$) if there
exists a positive constant $C>0$ such that for all finite
sequences $(x_n)$ of mutually disjoint elements in
$E$\hfil \break
\centerline { $( \sum \Vert x_n \Vert _{_{E}}^q )^{1/q}\leq
C\Vert \sum x_n\Vert _{_{E}} ,$}
\centerline {$ ({\rm resp.}  C( \sum \Vert x_n \Vert _{_{E}}^p
)^{1/p}\geq \Vert \sum x_n\Vert _{_{E}} ).$}

It is clear that if $E$ is $q-$concave for some $q\in
(0,\infty )$, then $E$ satisfies a lower $q-$estimate; if
$E$ satisfies a lower $q-$estimate for some $q\in (1,
\infty )$ then $E$ is $r-$concave, for all
$r>q$ [LT]. However, as shown in [Cr] Theorem
3.4(i), the Lorentz space $L_{q,p}(\jr^+ ),\ 1\leq p< q$
satisfies a lower $q-$ estimate, but is not $q-$concave.

We denote by \ $\nm$ \ a semi-finite von Neumann algebra in the Hilbert
space \
$H$ with given normal faithful semi-finite trace \ $\tau$. \ If \ $x$ \
is a
(densely defined) self-adjoint operator in \ $H$ \ and if \
$x=\int_{(-\infty, \infty)}sde^x_s$ \ is its spectral decomposition
then, for
any Borel subset \ $B\subseteq \jr$, \ we denote by \ $\ch_{_B}(x)$
\ the
corresponding spectral projection \ $\int_{(-\infty,
\infty)}\ch_{_B}(s)de^x_s$. \ A closed densely defined linear operator \
$x$ \
in \ $H$ \ affiliated with \ $\nm$ \ is said to be \ $\tau$-{\it
measurable} if
and only if there exists a number \ $s\ge0$ \ such that
$$\tau(\ch _{(s,\infty)}(\vert x\vert )) \ < \ \infty.$$
The set of all \ $\tau$-measurable operators will be denoted by \
$\nmt$. \ The
set \ $\nmt$ \ is a *-algebra with sum and product being the respective
closures of the  algebraic sum and product.  For \ $x\in\nmt$, \ the
{\it generalized singular value function} (or {\it decreasing
rearrangement}) \
$\mu .(x)$ \ of \ $x$ \ is defined by
$$\mu_t(x) \ = \ \inf\{ s\ge 0 \ : \ \tau(\ch_{(s,\infty)}(\vert x\vert
))\le
 t \}, \qquad t\ge 0.$$
It should be noted that the decreasing rearrangement $\mu (x)$ is the
right
continuous, non-increasing inverse of the distribution function $\lambda
_.(x)$, where $\lambda _{_{s}}(x)=\tau(\ch _{(s,\infty)}(\vert x\vert
)), s\geq 0.$  It
may be shown (cf. [FK]) that if $x\in \nmt$ and if $t>0$ then
$$\mu_t(x)=\inf
\{\Vert xe\Vert _{_{\nm }}:e\in \nm {\rm \  a\  projection\  with\ }
\tau
(1-e)\leq t \}.$$ It follows that \ $\mu(x)$ \ is a decreasing,
right-continuous
function on the half-line \ $\r=[0,\infty)$, and that  $\mu (x)=\mu
(x^*)=\mu
(|x|)$, for all $x\in \nmt $. For basic properties of decreasing
rearrangements
of measurable operators, we refer to [FK].  The topology defined by the
translation invariant metric \ $d$ \ on  \ $ \widetilde {\cal M}$ \
obtained by
setting  $$d(x,y)=\inf \lbrace t\ge 0:\mu _t(x-y)\leq t\rbrace,\qquad
{\rm for}\
x,y\in \widetilde {\cal M},$$ is called the {\it measure topology}. It
is shown
in  [Ne]
 and [Te] that  \ $ \widetilde {\cal M}$ \
equipped with the measure topology is a complete,
Hausdorff, topological *-algebra in which  \ ${\cal M}$ \ is
dense. It is not difficult to see that a net $(x_i)_{i\in I}$ in
$\widetilde
{\cal M}$ converges to  a measurable operator $x$ for the measure
topology if
and only if for every $\epsilon >0,\ \delta >0$ there exists $i_0\in I$
such
that whenever $i\geq i_0$ there exists a projection $e\in {\cal M}$ such
that
 $$ \Vert (x_i-x)e\Vert _{_{{\cal M}}}<\epsilon \quad {\rm and }\quad
\tau
(1-e)<\delta. $$ We denote by $\nmt _0$ the set of all $x\in \nmt $ such
that $\mu _t(x)\to 0$ as $t\to \infty .$

If \ $\nm$ \ is the space \ ${\cal L}(H)$ \ of all bounded linear
operators on \ $H$ \
and \ $\tau$ \ is the standard trace, then \ $\nmt={\cal L}(H)$ \ and \
$x\in
{\cal L}(H)$ \ is compact if and only if \ $\mu_t(x)\to 0$ \ as \ $t\to
\infty$,
\ in which case for each \ $n=0, 1, 2,\dots$, $$\mu_n(x) \ = \ \mu_t(x),
\
\qquad t\in [n, n+1),$$ and \ $\{ \mu_n(x)\} _{n=0}^\infty$ \ is the
usual
singular value sequence of \ $x$ \ in decreasing order, counted
according to multiplicity [GK].

We identify the space \ $L^\infty(\r)$ \ of all bounded complex-valued
Lebesgue
measurable functions on the half-line \ $\r$ \ as a commutative von
Neumann
algebra acting by multiplication on the Hilbert space \ $L^2(\r)$ \ with
trace
given by integration with respect to Lebesgue measure. In this case the
\ $\tau$-measurable operators coincide with those complex measurable
functions
\ $f$ \ on \ $\r$ \ which are bounded except on a set of finite measure.
In
this example, the generalized singular value function, which we continue
to
denote by \ $\mu (f)$, \ coincides with the familiar right-continuous
decreasing
rearrangement of the function \ $|f|$. \ See, for example, [KPS].

Following [X], by a symmetric quasi-Banach function space on the
positive
half-line $\r $ is meant a quasi-Banach lattice $E$ of measurable
functions
with the following properties: (i) $E$ is an order-ideal in the linear
space
$L^0(\r)$ of all finite almost-everywhere measurable functions on $\r$;
(ii) $E$ is rearrangement-invariant in the sense that  whenever $x\in E$
and
$y$ is a
measurable function with $\mu (y)=\mu (x)$, it follows that $y\in E $
and
$\Vert y\Vert _{_{E}}= \Vert x\Vert _{_{E}};$ (iii) $E$ contains all
finitely-supported simple functions.  The norm on $E$ is said to be a
$\sigma -${\it Fatou} norm if  $$ 0\leq x_n\uparrow x,\quad x_n,x\in E \
{\rm
implies }\  \Vert x_n\Vert _{_{E}} \uparrow \Vert x\Vert _{_{E}}.$$
The norm on $E$ is said to be {\it order-continuous} if
$$ 0\leq x_{\beta }\downarrow _{_{\beta }} 0\ {\rm implies }\ \Vert
x_{_{\beta }}\Vert _{_{E}}\downarrow _{_{\beta }} 0.$$
We remark that if $E$ has order-continuous norm, then $\lim _{_{t\to
\infty
}}\mu _{_{t}}(x) =0 $ for all $x\in E.$ It is clear that if $E$ has
order-continuous norm, then the norm on $E$ is a $\sigma -$Fatou norm.
Further,
it is not difficult to see that if $E$ is $q$-concave for some $0<q<\infty
$, then
the norm on $E$ is order-continuous.  Now suppose that the symmetric
quasi-Banach
function space $E$ is $\alpha-$convex for some $0<\alpha <\infty $ and
that
$M^{(\alpha)}(E)=1.$  Let $F$ be the Banach lattice $E^{(1/\alpha )}$.
If the
norm on $E$ is a $\sigma-$Fatou norm, then so is the norm on $F$ and
consequently the natural embedding of $F$ into the second associate
space ( K\"othe bidual)
$F^{\prime \prime}$ is an isometry,
 by a well known result of Lorentz and Luxemburg [LT] 1.b.18.
Consequently, by [KPS] II.4.6, there exists a family $W$ of non-increasing
functions on
$\r$ such that  $$ \Vert x \Vert _{_{F}}=\sup \left\{ \int _{[0, \infty )}\mu
_t
(x)w(t)dt\ :\ w\in W \right\}\eqno (1.1)$$ for all $x\in F.$  Consequently,
$$\Vert x \Vert _{_{E}}^{\alpha }=\sup \left\{ \int _{[0, \infty )}\mu
_t^{\alpha }
(x)w(t)dt\ :\ w\in W \right\}\eqno (1.2)$$
for all $x\in E.$

Unless stated otherwise, $E$ will always denote a symmetric
quasi-Banach function space on the positive half-line $\r $ and
$(\nm,\tau)$ a
semifinite von Neumann algebra on a Hilbert space $H$, with 
faithful normal semi-finite trace $\tau $. We define the symmetric space
$E(\nm)$
of measurable operators associated with $E$ by setting
$$E(\nm )=\{x\in \nmt\ : \mu (x)\in E \}$$ and $$\Vert x\Vert
_{_{E(\nm)}}=\Vert \mu (x)\Vert _{_{E}},\quad x\in E(\nm ).$$
It follows immediately that $\Vert x\Vert _{_{E(\nm)}}=
\Vert \ |x|\ \Vert _{_{E(\nm)}}=\Vert x^*\Vert _{_{E(\nm)}}$ for all
$x\in E(\nm) .$
 It is shown in [X] Lemma 4.1 that if $E$ has $\sigma-$Fatou norm and
is
$\alpha-$convex  for some $0<\alpha <\infty$ with $M^{(\alpha )}(E)=1$,
then the
functional $\Vert \cdot \Vert _{_{E(\nm)}}$ is a norm if $\alpha \geq 1$
and
an $\alpha -$norm if $0<\alpha <1.$ This means that if $x,y \in E(\nm
)$, then
$$\Vert x+y\Vert _{_{E(\nm)}}^{\alpha }\leq \Vert x\Vert _{_{E(\nm)}}^
{\alpha }+\Vert y\Vert _{_{E(\nm)}}^ {\alpha }.$$
Equipped with this norm or $\alpha -$norm, the space $E(\nm)$ is
complete. For related results, see [Ov], [DDP1,2].

If \ $x,y\in \widetilde {\cal M}$, \  we say that \ $x$ \ is {\it
submajorized
by} \ $y$, \  written \  $x\prec \prec y$, \  if and only if
$$\int _0 ^\alpha \mu _t(x)dt\le \int_0^\alpha
\mu_t(y)dt,\qquad {\rm for\ all}\ \alpha \ge 0.$$
If $E$ is a symmetric Banach function space with $\sigma-$Fatou norm and
if
$x,y\in E$ satisfy $x\prec \prec y$, then $\Vert x\Vert _{_{E}}\leq
\Vert y\Vert _{_{E}}$; if $E$
has order-continuous norm or is {\it maximal} [KPS] in the sense that
the natural
embedding of $E$ into its second associate space is an isometric
surjection,
then $E$ has the following property: if $y\in E$ and if $x\prec \prec
y$, then
$x\in E $ and $\Vert x\Vert _{_{E}}\leq \Vert y\Vert _{_{E}}.$
Suppose now that $E$ is $\alpha -$convex with $M^{(\alpha )}(E)=1$ for some
$\alpha
\in (0, \infty),$ and that $0<q<\infty $. If $x,y\in \nmt,$ it follows from
[FK],
Theorem 4.2(iii) that  $$ \mu ^{\alpha q} (xy)\prec \prec \mu ^{\alpha
q}(x)\mu
^{\alpha q }(y).$$ It now follows from the representation (1.2), the
usual
H\"older inequality and a well-known rearrangement inequality of Hardy
[KPS] II.2.18,
that for $0<q_{_{0}},q_{_{1}}, q<\infty $
with $1/q=1/q_{_{0}}+1/q_{_{1}}$
$$\Vert xy\Vert _{_{E^{(q)}(\nm)}}\leq \Vert x\Vert
_{_{E^{(q_{_{0}})}(\nm)}}
\Vert y\Vert _{_{E^{(q_{_{1}})}(\nm)}}.\eqno (1.3)$$
We shall need the following result, which is well known for the case
that $E$
is a Banach lattice [DDP2], [CS].  Here, order-continuity of the norm on
$E(\nm)$
is defined in the obvious way.

\noindent {\bf Proposition 1.1}\quad {\sl Let $E$ be $\alpha
-$convex
with $M^{(\alpha )}(E)=1$ for some $\alpha \in (0, 1]$.
\item {\rm (i).}  If $E$ has order-continuous norm, then so does
$E(\nm)$.
\item {\rm(ii).} If $E$ has order-continuous norm, and if
$\{e_{_{\beta}}\}\subseteq \nm $ is any family of (self-adjoint)
projections for
which $e_{_{\beta}}\downarrow _{_{\beta }} 0$ then $\lim _{_{\beta
}}\Vert
xe_{_{\beta}}\Vert _{_{E(\nm )}}=0$ for every} $x\in E(\nm).$

\noindent {\bf Proof.}\quad  (i). Suppose that $E$ has order-continuous
norm
and that $0\leq x_{_{\beta}}\downarrow _{_{\beta }} 0$ holds in
$E(\nm)$.
As noted earlier, order-continuity of the norm on $E$ implies that $\lim
_{_{t\to \infty }}\mu _{_{t}}(x) =0 $ for all $x\in E.$  It follows from
[DDP2] (see also [CS]) that $ \mu (x_{_{\beta}})\downarrow _{_{\beta
}}0$ holds in $E$ and (i) now follows. \hfil\break
\noindent (ii). Suppose that $\{e_{_{\beta}}\}\subseteq \nm $ is any
family of
projections in $\nm$ with $e_{_{\beta }}\downarrow _{_{\beta }} 0$ and
let
$0\leq x\in \nm .$ Using (1.3) above, it follows that
$$\eqalign {\Vert xe_{_{\beta }}\Vert  _{_{E(\nm)}} &= \Vert
e_{_{\beta }} x^{1/2} x^{1/2}\Vert  _{_{E(\nm)}}  \cr
&\leq \Vert (x^{1/2})^2 \Vert  _{_{E(\nm)}} ^{1/2}
\Vert |e_{_{\beta }}x^{1/2}|^2 \Vert  _{_{E(\nm)}} ^{1/2} \cr
&=\Vert x \Vert  _{_{E(\nm)}} ^{1/2}\Vert x^{1/2}e_{_{\beta }}x^{1/2}
\Vert  _{_{E(\nm)}} ^{1/2} .\cr}$$
Now $e_{_{\beta }}\downarrow _{_{\beta }} 0$ in $\nm$ implies that
$ x^{1/2}e_{_{\beta }}x^{1/2} \downarrow _{_{\beta }} 0$ in $E(\nm)$
and the assertion of (ii) now follows from that of (i).

\noindent {\bf 2. A uniform Kadec-Klee property}
\smallskip

It is now convenient to localize the notion of convergence in measure as
follows.  If $\epsilon ,\delta >0$
and if $e$ is a projection in $\nm $
 with $\tau (e)<\infty $, then the family of all sets $N_{_{\epsilon
,\delta,e }}$ consisting of all $x\in \nmt $ such that $\mu
_{_{\delta}}(exe)<\epsilon $ form a neighbourhood base at $0$ for a
Hausdorff
linear topology on $\nmt $. This topology (cf.[ DDP2], Proposition 5.13)
will
be called the topology of {\it local convergence in measure} ($lcm$).
 If we observe that
$$(|exe|-\lambda 1)^+=(|exe|-\lambda e)^+$$
for all $\lambda >0,\ x\in \nmt$ and projections $e\in \nm$, then it
follows
that $\mu (exe)=\mu ^e(exe)$ for all $x\in \nmt $, where $\mu ^e$
denotes
decreasing rearrangement calculated relative to the von Neumann algebra
$e\nm
e$ with respect to the trace $\tau (e\cdot e).$  This remark follows
also as a
special case of  [Fa], Proposition 1.5(i). It follows that convergence
locally
in measure coincides with convergence for the measure topology relative
to
 $\left (e\nm e, \tau (e\cdot e)\right )$, for each projection $e\in \nm
$ with
$\tau (e)<\infty $.

We remark immediately that if $\nm $ is commutative, and identified with
the
von Neumann algebra of all multiplication operators given by bounded
measurable
functions on some localisable measure space then the preceding notion of
local convergence in measure reduces to the more familiar notion of
convergence
in measure on sets of finite measure in the underlying measure space. On
the
other hand, if \ $\nm$ \ is the space \ ${\cal L}(H)$ \ of all bounded
linear
operators on \ $H$ \ and \ $\tau$ \ is the standard trace, then
convergence
locally in measure is precisely convergence for the weak operator
topology.

\noindent {\bf Proposition 2.1 }\quad {\sl Assume that E is
$\alpha-convex $ for some $0<\alpha \leq 1$ with $M^{(\alpha )}(E)=1$.  If
$E$ has
order-continuous norm, and if $\{x_{_{n}}\}$ is a sequence in $E(\nm)$
which
converges to $x\in E(\nm)$ locally in measure, then for every $\epsilon
>0$,
there exists a projection $e$ in $\nm$ with $\tau (e) <\infty $ and a
subsequence
$\{y_{_{n}}\} \subseteq \{x_{_{n}}\}$ such that}
$$ \Vert e(y_{_{n}}-x)e\Vert _{_{E(\nm)}}\to _n 0\ {\sl and}\ {\rm max }
\{\Vert x(1-e)\Vert _{_{E(\nm)}},\Vert x^*(1-e)\Vert
_{_{E(\nm)}}\}<\epsilon
.$$
\noindent {\bf Proof}\quad The semi-finiteness of $\nm$ implies that the
family $\{ e_{_{\beta }} \}$ of all projections in $\nm$ of finite trace
satisfies
$0\leq e_{_{\beta}}\uparrow_{_{\beta }} 1$. Since  $E$ has
order-continuous
norm, it follows from Proposition 1.1 that
$${\rm max}\{\Vert x(1-e_{_{\beta }})\Vert _{_{E(\nm)}},\Vert
x^*(1-e_{_{\beta }})\Vert _{_{E(\nm)}}\}\to _\beta 0.$$
Consequently, there exists a projection $e_{_{0}}\in \nm $ with finite
trace
such that
$${\rm max}\{\Vert x(1-e_{_{0 }})\Vert _{_{E(\nm)}}^{\alpha },\Vert
x^*(1-e_{_{0}})\Vert _{_{E(\nm)}}^{\alpha }\}<\epsilon ^{\alpha } /2.$$
We next observe that by passing to a subsequence and relabelling if
necessary,
there exists a sequence $e_{_{m}}$ of projections in $\nm$ with
$e_{_{m}}\uparrow _{_{m}}e_{_{0}}$ such that,for all $m=1,2,\dots $ ,
$$ \Vert e _{_{0}}(x_{_{n}}-x)e_{_{m}}e_{_{0}}\Vert _{_{\infty }}\to
_n0,$$
where $\Vert \cdot \Vert _{_{\infty }}$ denotes the usual operator norm
on $\nm
$.  While the proof of this assertion is a standard argument, we include
the
details for the sake of completeness. In fact, since $
e_{_{0}}(x_{_{n}}-x)e_{_{0}}\to 0$ for the measure topology given by
$\left (e_{_{0}}\nm e_{_{0}},\tau( e_{_{0}}\cdot e_{_{0}})\right )$, it
follows, by passing to a subsequence if necessary, that there exist
projections
$e_{_{k}}^{\prime }\leq e_{_{0}},k=1,2,\dots $ such that
$$\tau(e_{_{0}}-e_{_{k}}^{\prime } ) <2^{-k}\quad {\rm and }\quad \Vert
e_{_{0}}(x_{_{k}}-x)e_{_{k}}^{\prime }e_{_{0}}\Vert _{_{\infty
}}<2^{-k}.$$
For $m=1,2,\dots $, define $e _{_{m}}=\inf _{k\geq m}e_{_{k}}^{\prime
}.$
Observe that
$$e_{_{0}}-e_{_{m}}=e_{_{0}}-\inf _{k\geq m}e_{_{k}}^{\prime }=\sup
_{k\geq
m}(e_{_{0}}-e_{_{k}}^{\prime }), $$
and consequently
$$\tau (e_{_{0}}-e_{_{m}})\leq \sum _{k=m}^{\infty }\tau
(e_{_{0}}-e_{_{k}}^{\prime })\leq 2^{-m+1}, $$
for all $m=1,2,\dots $. This implies that $ e_{_{m}}\uparrow _{_{m}}e
_{_{0}}.$
Finally, if $n\geq m$, then
$$ \Vert e_{_{0}}(x_{_{n}}-x)e_{_{m}}e_{_{0}}\Vert _{_{\infty }}=
\Vert e_{_{0}}(x_{_{n}}-x)e_{_{n}}^{\prime }e_{_{m}}e_{_{0}}\Vert
_{_{\infty
}}\leq \Vert e_{_{0}}(x_{_{n}}-x)e_{_{n}}^{\prime }e_{_{0}}\Vert
_{_{\infty }}
\leq 2^{-n},$$
and the assertion follows.  Order-continuity of the norm on $E$ again
implies
that there exists a natural number $r$ such that
$${\rm max}\{\Vert x(e_{_{0 }}-e_{_{m}})\Vert _{_{E(\nm)}}^ {\alpha
},\Vert
x^*(e_{_{0}}-e_{_{m}})\Vert _{_{E(\nm)}}^{\alpha }\}<\epsilon ^{\alpha
}/2,\quad
m\geq r.$$ Consequently,
$$ \Vert x(1-e_{_{r}})\Vert _{_{E(\nm )}}^{\alpha }\leq
\Vert x(1-e_{_{0}})\Vert _{_{E(\nm )}}^{\alpha }+ \Vert
x(e_{_{0}}-e_{_{r}})
\Vert _{_{E(\nm)}}^{\alpha }<\epsilon ^{\alpha } ,$$
and similarly $\Vert x^*(1-e_{_{r}})\Vert _{_{E(\nm )}}< \epsilon .$ We
now
observe that
$$\Vert e_{_{r}}(x_{_{n}}-x)e_{_{r}} \Vert _{_{\infty }}=
\Vert e_{_{0}} e_{_{r}}(x_{_{n}}-x)e_{_{r}}e_{_{0}} \Vert _{_{\infty }}
\leq \Vert e_{_{0}}(x_{_{n}}-x)e_{_{r}}e_{_{0}} \Vert _{_{\infty }} $$
and so
$$ \Vert e_{_{r}}(x_{_{n}}-x)e_{_{r}} \Vert _{_{E(\nm)}}\leq
\Vert e_{_{r}}(x_{_{n}}-x)e_{_{r}} \Vert_{_{\infty }} \Vert
e_{_{0}}\Vert _{_{
E(\nm)}} $$
and this suffices to complete the proof of the proposition.

\noindent {\bf Lemma 2.2}\quad {\sl Let $0\leq x\in \nmt $, let $e$ be a
projection in $\nm $ with $\tau (e) <\infty $ and let $f=1-e$. If
$0<\gamma
<\infty $, then for all $0\leq a_{_{1}},a_{_{2}} \in \jr$},
$$\int _0^{a_{_{1}}}\mu _{_{t}}^{\gamma }(exe)dt+
\int _0^{a_{_{2}}}\mu _{_{t}}^{\gamma }(fxf)dt\leq
\int _0^{a_{_{1}}+a_{_{2}}}\mu _{_{t}}^{\gamma }(exe+fxf)dt. $$
\noindent{\bf Proof }\quad We begin by observing that, for all $s>0$,
$$\ch _{(s,\infty)}(\vert exe+fxf \vert )=
\ch _{(s,\infty)}(\vert exe\vert )+\ch _{(s,\infty)}(\vert fxf\vert )$$
and so
$$\lambda _{_{s}}(exe+fxf)=\lambda _{_{s}}(exe)+\lambda _{_{s}}(fxf)$$
for all $s>0$.  Using the fact that the functions $\lambda , \mu $ are
essential inverses, it follows that
$$ \int _0^u\mu _{_{t}}^{\gamma}(z)dt=\int _0^{\infty } u\wedge \lambda
_{_{s^{1\over
\gamma }}}(z)ds $$
 for all $z\in \nmt $, and for all
$0\leq u\in \jr$.  Consequently,
$$\eqalign {\int _0^{a_{_{1}}}\mu _{_{t}}^{\gamma }(exe)dt+
\int _0^{a_{_{2}}}\mu _{_{t}}^{\gamma }(fxf)dt&=
\int _0^{\infty }\left ( a_{_{1}}\wedge \lambda _{_{s^{1\over \gamma
}}}(exe)+
a_{_{2}}\wedge \lambda _{_{s^{1\over \gamma }}}(fxf) \right )ds\cr
&\leq \int _0^{\infty  }(a_{_{1}}+a_{_{2}})\wedge \left (\lambda
_{_{s^{1\over
\gamma}}}(exe)+  \lambda _{_{s^{1\over \gamma }}}(fxf)\right ) ds\cr
&= \int _0^{\infty  }(a_{_{1}}+a_{_{2}})\wedge \lambda _{_{s^{1\over
\gamma}}}(exe+fxf) ds\cr
&=\int _0^{a_{_{1}}+a_{_{2}}}\mu _{_{t}}^{\gamma }(exe+fxf)dt.\cr} $$

We remark that the preceding proof is a variant of the argument of [Fr]
Lemma 16.

\noindent {\bf Lemma 2.3}\quad {\sl Let $x\in \nmt $, let $e$ be a
projection in $\nm $ with $\tau (e) <\infty $ and let $f=1-e.$ If
$0<\gamma
<\infty $, then}
$$\mu ^{\gamma}(e\vert x\vert e)+\mu ^{\gamma }_{_{(\cdot )-\tau
(e)}}(f\vert
x\vert f)\prec \prec \mu ^{\gamma }(e\vert x\vert e+f\vert x\vert f).$$
\noindent {\bf Proof}\quad It suffices to show that if $X$ is a Lebesgue
measurable subset of $\r $, then
$$\int _X \left [\mu _{_{t}}^{\gamma }(e\vert x\vert e)+\mu _{_{t-\tau
(e)}}^{\gamma}(f\vert x\vert f)\right ]dt\leq \int_0^{\vert X\vert }\mu _{_{t}}
 ^{\gamma}(e\vert x\vert e+f\vert x\vert f)dt, $$
where $ \vert X\vert $ denotes the Lebesgue measure of $X$. Now observe
that
$$\eqalign {\int _X \left [\mu _{_{t}}^{\gamma }(e\vert x\vert e)+\mu
_{_{t-\tau
(e)}}^{\gamma}(f\vert x\vert f)\right ]dt
&=\int _{_{X\cap [0,\tau (e)]}}\mu _{_{t}}^{\gamma }(e\vert x\vert e)dt
+
\int _{_{X\cap [\tau (e),\infty )}}\mu _{_{t-\tau (e)}}^{\gamma }(f\vert
x\vert
f)dt\cr
&\leq \int _{_{0}}^{\vert X\cap [0,\tau (e)]\vert }\mu _{_{t}}^{\gamma
}(e\vert
x\vert e)dt + \int _{_{0}}^{\vert X\cap [\tau (e),\infty )\vert }
\mu _{_{t}}^{\gamma }(f\vert x\vert f)dt.\cr }$$
By observing that
$$\vert X\cap [0,\tau (e)]\vert +\vert X\cap [\tau (e),\infty )\vert
=\vert
X\vert $$
the desired assertion now follows from Lemma 2.2 preceding.

We shall need the following submajorization inequality, which is
essentially
well known [Si], [GK]. For sake of completeness, we include a proof
which is a
simple modification of an argument of J. Arazy [Si].

\noindent {\bf Lemma 2.4}\quad {\sl Let $0\leq x\in \nmt $, let $e$ be a
projection in $\nm $ and let $f=1-e.$  If $1\leq \gamma <\infty $, then}
$$\mu ^{\gamma }(exe+fxf)\prec \prec \mu ^{\gamma }(x).$$
\noindent {\bf Proof } Set $v=e-f$; note that $v\in \nm$ is unitary, and
that
$$
exe+fxf={1\over 2}(x+vxv^*).$$It follows that
$$ \mu (exe+fxf)=\mu ( {1\over 2}(x+vxv^*))\prec \prec {1\over 2}(\mu
(x)+
\mu (vxv^*))=\mu(x).$$
This shows that the assertion of the Lemma is true if $\gamma =1$. The
case
 where $\gamma >1 $ follows from the case that $\gamma =1$, and the fact
that the
function $(\cdot )^{\gamma } $ is convex if $\gamma >1$.

It is not difficult to see that the assertion of Lemma 2.5 fails if
$0<\gamma< 1 $,  by taking  $\nm={\cal L}({\jc}^2)$ ,
$$ x=\pmatrix{
5&4\cr
4&5\cr},$$
and by letting $e,f$ be the natural coordinate projections. We omit the
details.

We now come to the principal estimate needed for our main result. The
inequality given below is a refinement of a similar estimate
for trace class operators given by Arazy [Ar] and Lennard [Le1].

\noindent {\bf Proposition 2.5}\quad {\sl Let $E$ be $\alpha -$convex
with
$M^{(\alpha )}(E)=1$ for some $0<\alpha \leq 1$ and suppose that $E$
satisfies a
lower $q$-estimate with constant 1 for some finite $q\geq \alpha $. If
$k=2q/\alpha $, then for all $y\in E (\nm)$ and all projections $e,f\in
\nm
$ with $e+f=1$ and $\tau (e)<\infty $, it follows that}
$$\Vert eye\Vert ^k_{_{E(\nm)}}+\Vert eyf\Vert ^k_{_{E(\nm)}}+
\Vert fye\Vert ^k_{_{E(\nm)}}+\Vert fyf\Vert ^k_{_{E(\nm)}}\leq
\Vert y\Vert ^k_{_{E(\nm)}}.$$
\noindent {\bf Proof}\quad The stated assumptions on $E$ imply that
$E^{1/\alpha }$ is a Banach lattice which satisfies a lower $(q/\alpha
)-$estimate, with constant $1$.  Let $0\leq w\in E^{(1/\alpha )}(\nm)$
and let
$e,f$ be as stated. Using the lower estimate, it follows that
$$\eqalign {\Vert ewe\Vert ^{q/\alpha }_{_{E^{(1/\alpha )}(\nm)}}+
\Vert fwf\Vert ^{q/\alpha }_{_{E^{(1/\alpha )}(\nm)}}
&=\Vert \mu (ewe)\Vert ^{q/\alpha }_{_{E^{(1/\alpha )}}}+
\Vert \mu _{_{(\cdot ) -\tau(e)}}(fwf)\Vert ^{q/\alpha }_{_{E^{(1/\alpha
)}}}\cr
&\leq \Vert \mu (ewe)+\mu _{_{(\cdot )-\tau (e)}}(fwf)\Vert ^{q/\alpha
}_{_{E^{(1/\alpha )}}}\cr
&\leq\Vert \mu (ewe+fwf)\Vert ^{q/\alpha }_{_{E^{(1/\alpha )}}}\cr
&\leq \Vert \mu (w)\Vert ^{q/\alpha }_{_{E^{(1/\alpha )}}}=
\Vert w\Vert ^{q/\alpha }_{_{E^{(1/\alpha )}(\nm)}},\cr }$$
where the last two inequalities follow from successive
application of Lemmas 2.3, 2.4. It now follows that for all $z\in
E^{(2/\alpha
)}(\nm)$ that
$$\Vert ez\Vert ^k _{_{E^{(2/\alpha )}(\nm )}}+\Vert fz\Vert ^k
_{_{E^{(2/\alpha )}(\nm )}}
\leq\Vert z\Vert ^k _{_{E^{(2/\alpha )}(\nm )}}.\eqno (*)$$
In fact,
$$\eqalign {\Vert ez\Vert ^k _{_{E^{(2/\alpha )}(\nm)}}+\Vert fz\Vert ^k
_{_{E^{(2/\alpha )}(\nm)}}
&=\Vert (ezz^*e)^{1/2}\Vert ^k _{_{E^{(2/\alpha )}(\nm)}}+\Vert
(fzz^*f)^{1/2}
\Vert ^k_{_{E^{(2/\alpha )}(\nm)}}\cr
&=\Vert ezz^*e\Vert ^{k/2} _{_{E^{(1/\alpha )}(\nm)}}+\Vert fzz^*f
\Vert ^{k/2}_{_{E^{(1/\alpha )}(\nm)}}\cr
&\leq \Vert zz^* \Vert^{k/2}_{_{E^{(1/\alpha )}(\nm )}}=
\Vert z \Vert^k_{_{E^{(2/\alpha )}(\nm )}}.\cr }$$
Set $\beta=2/\alpha$ and choose $\gamma $, $1<\gamma \leq \infty $, such that
$2/\beta +1/\gamma =1$, which is possible since $\beta \geq 2$. Suppose now that
$y\in E(\nm)$ and
set
$$ a=u\vert y\vert ^{1/\beta },\quad b=\vert y\vert ^{1/\gamma },\quad
c=\vert y\vert ^{1/\beta }, $$
where $y=u\vert y\vert $ is the polar decomposition of $y$. It follows
that $a,c\in E^{(\beta )}(\nm)$ and $b\in  E^{(\gamma )}(\nm)$ (if
$\gamma=\infty$, we take  $b$ to $1$ and $E^{(\gamma )}(\nm)$ to be
$\nm$.
Note that
$$\Vert a\Vert^{\beta }_{_{E^{(\beta )}(\nm)}}=
\Vert b\Vert^{\gamma }_{_{E^{(\gamma )}(\nm)}}=
\Vert c\Vert^{\beta }_{_{E^{(\beta )}(\nm)}}=
\Vert y\Vert_{_{E(\nm)}}.\eqno (**)$$
We apply the H\"older-type inequality (1.3) to obtain
$$\eqalign {\Vert eye\Vert _{_{E(\nm)}}
&=\Vert eabce\Vert _{_{E(\nm)}}\cr
&\leq \Vert ea\Vert _{_{E^{(\beta )}(\nm)}}
\Vert b\Vert _{_{E^{(\gamma )}(\nm)}}
\Vert ce\Vert _{_{E^{(\beta )}(\nm)}}.\cr}$$
Similarly,
$$\eqalign {\Vert eyf\Vert _{_{ E(\nm)}}
&\leq \Vert ea\Vert _{_{ E^{(\beta )}(\nm)}}
\Vert b\Vert _{_{ E^{(\gamma )}(\nm)}}
\Vert cf\Vert _{_{ E^{(\beta )}(\nm)}},\cr
\Vert fye\Vert _{_{ E(\nm)}}
&\leq \Vert fa\Vert _{_{ E^{(\beta )}(\nm)}}
\Vert b\Vert _{_{ E^{(\gamma )}(\nm)}}
\Vert ce\Vert _{_{E^{(\beta )}(\nm)}},\cr
\Vert fyf\Vert _{_{ E(\nm)}}
&\leq \Vert fa\Vert _{_{ E^{(\beta )}(\nm)}}
\Vert b\Vert _{_{ E^{(\gamma )}(\nm )}}
\Vert cf\Vert _{_{ E^{(\beta )}(\nm )}}.\cr} $$
Consequently, using first (*) and then (**) we obtain that
$$\eqalign {&\Vert eye\Vert ^k_{_{E(\nm)}}+\Vert eyf\Vert
^k_{_{E(\nm)}}+
\Vert fye\Vert ^k_{_{E(\nm)}}+\Vert fyf\Vert ^k_{_{E(\nm)}}\cr
&\leq \left (\Vert ea\Vert ^k_{_{E^{(\beta )}(\nm)}}+
\Vert fa\Vert ^k_{_{E^{(\beta )}(\nm)}}\right )
\Vert b\Vert ^k_{_{E^{(\gamma )}(\nm)}}
\left (\Vert cf\Vert ^k_{_{E^{(\beta )}(\nm)}}+
\Vert ce \Vert ^k_{_{E^{(\beta )}(\nm)}}\right )\cr
&\leq \Vert a\Vert ^k_{_{E^{(\beta )}(\nm)}}
\Vert b\Vert ^k_{_{E^{(\gamma )}(\nm)}}
\Vert c\Vert ^k_{_{E^{(\beta )}(\nm)}}\cr
&=\Vert y\Vert ^{k(2/\beta+1/\gamma )}_{_{E(\nm)}}\cr
&=\Vert y\Vert ^k_{_{E(\nm)}},\cr }$$
and the proof of the proposition is complete.

Before 
proceeding, it is convenient to
introduce the
following terminology.

\noindent {\bf Definition 2.6}\quad {\sl Let $\left (X,\Vert\cdot \Vert
_{_{X}}\right )$ be an $\alpha-$normed quasi-Banach space, for some $0<
\alpha \leq 1$, with unit ball ${\bf B}_{_{X}}$, and let ${\cal T}$ be a
topological vector space topology on $X$ that is weaker than the
quasi-norm
topology. The space  $(X,\Vert\cdot \Vert _{_{X}})$
 is said to have the uniform Kadec-Klee-Huff property with respect to
${\cal T}$ ( denoted $UKKH({\cal T})$)
 if and only if for every $\epsilon >0$ there
exists $\delta=\delta(\epsilon ) \in (0,1)$
 such that whenever $\{x_{_{n}}\}$ is a
sequence
in ${\bf B}_{_{X}}$ with $x_{_{n}}\to _n x \in X$ with respect to ${\cal
T}$ and
$\Vert x \Vert _{_{X}}>1-\delta $, it follows that for some
 positive integer $N$,}
$$\sup _{n, m\geq N}\Vert x_{_{n}}-x_{_{m}} \Vert _{_{X}}\leq \epsilon
.$$

Let us immediately note that the conclusion in the statement of the
preceding
definition is stronger than the usual formulation of this
property given by several authors (cf [H],
 [Le1,2], [DL1,2]) in which $\sup $ is
replaced by $\inf $. However, a  minor modification of the argument of
[BDDL]
Proposition 1.2 shows that our (apparently) stronger formulation is, in
fact,
equivalent to that already in the literature.

 We may now state the main result of this paper.

\noindent {\bf Theorem 2.7}\quad {\sl Let $E$ be a symmetric quasi-Banach
functions space on the positive half-line $\r $ If $E$ is $\alpha -$convex with
$M^{(\alpha )}(E)=1$ for some $0<\alpha \leq 1$ and if $E$ satisfies a lower-$q$
estimate with constant 1 for some finite $q\geq \alpha $, then $\left (E(\nm),
\Vert \cdot \Vert _{_{E(\nm)}}\right )$ has the uniform Kadec-Klee-Huff property
with respect to local convergence in measure.}

\noindent {\bf Proof }\quad Suppose that $\delta \in (0,1) $ and
that $\{x_{_{n}}\}^{\infty }_{_{n=1}}$is a sequence in the unit ball
of $E(\nm)$ satisfying $x_{_{n}}\to _n x$ locally in measure for some
$x\in
E(\nm)$ with $\Vert x\Vert_{_{E(\nm)}}>1-\delta .$ By [BDDL],
Proposition 1.2,
it suffices to show that
 for some continuous, invertible function $\epsilon = \epsilon(\delta)$ mapping
 $(0,\delta_0)$ into $(0,\infty)$, where $\delta_0 > 0$ and
 $\epsilon(\delta) \to 0$ as $\delta \to 0$, we have
$$\inf _{n\neq m}\Vert x_{_{n}}-x_{_{m}}\Vert_{_{E(\nm)}}<\epsilon .$$
 By passing to a subsequence and relabelling, if necessary, and
appealing to
Proposition 2.1, we may assume that there exists a projection
$e\in \nm $ with $\tau (e)<\infty $ such that, setting $f=1-e$,
$$\Vert e(x_{_{n}}-x)e\Vert _{_{E(\nm)}}\to _n 0 \eqno(i)$$ and
$${\rm max}\left \{ \Vert xf\Vert _{_{E(\nm)}}, \Vert x^*f\Vert
_{_{E(\nm)}}\right \}<\left [{1\over 3}\left (\Vert x\Vert
_{_{E(\nm)}}^{\alpha }-(1-\delta )^{\alpha }\right )\right ]^
{1/\alpha }.\eqno (ii)$$
Since $\Vert \cdot \Vert _{_{E(\nm)}}$ is an $\alpha -$norm, we
obtain that
$$\Vert x_{_{n}}-x_{_{m}}\Vert ^{\alpha }_{_{E(\nm)}}\leq
\Vert e(x_{_{n}}-x_{_{m}})e\Vert ^{\alpha }_{_{E(\nm)}} +R(x_{_{n}})
+R(x_{_{m}}),$$
where
$$R(y)=\Vert eyf\Vert ^{\alpha }_{_{E(\nm)}}+
\Vert fye\Vert ^{\alpha }_{_{E(\nm)}}+
\Vert fyf\Vert ^{\alpha }_{_{E(\nm)}},$$
for all $y\in E (\nm )$.  Setting $k=2q/\alpha , s=k/\alpha $ and $s^{\prime
}=s/(s-1), $ and applying Proposition 2.5, we obtain, for each $y\in \nm
$,
$$\eqalign {R(y)&\leq  3^{1/s^\prime }\left (\Vert eyf\Vert
^k_{_{E(\nm)}}+
\Vert fye\Vert ^k_{_{E(\nm)}}+\Vert fyf\Vert ^k_{_{E(\nm)}}\right
)^{1/s}\cr
&\leq  3^{1/s^\prime }\left (\Vert y\Vert ^k_{_{E(\nm)}} -
\Vert eye\Vert^k_{_{E(\nm)}}\right )^{1/s}.\cr }$$Using (ii), it now
follows
that
$$ \eqalign {\Vert exe\Vert ^{\alpha }_{_{E(\nm)}}
&\geq \Vert x\Vert ^{\alpha }_{_{E(\nm)}}-\Vert exf\Vert ^{\alpha
}_{_{E(\nm)}}-
\Vert fxe\Vert ^{\alpha }_{_{E(\nm)}}-
\Vert fxf\Vert ^{\alpha }_{_{E(\nm)}}\cr
&\geq \Vert x\Vert ^{\alpha }_{_{E(\nm)}}
-3{\rm max}\left \{ \Vert xf\Vert _{_{E(\nm)}}^ {\alpha }, \Vert
x^*f\Vert
_{_{E(\nm)}}^{\alpha }\right \}\cr
&>\Vert x\Vert ^{\alpha }_{_{E(\nm)}}-
\left (\Vert x\Vert ^{\alpha }_{_{E(\nm)}}-(1-\delta)^{\alpha }\right
)\cr
&=(1-\delta)^{\alpha },\cr }$$
so that $ \Vert e x e \Vert _{_{E(\nm)}}> 1-\delta .$  From (i) and the
uniform continuity of the quasi-norm $\Vert \cdot \Vert _{_{E(\nm)}}$
 on the closed unit ball of $E(\nm )$, it follows that
$$\Vert ex_{_{n}}e\Vert _{_{E(\nm)}}\to _n
\Vert exe\Vert _{_{E(\nm)}} \ ;$$
 and so there exists a natural number $n_{_{0}}$ such that
$$\Vert ex_{_{n}}e\Vert _{_{E(\nm)}}>1-\delta,\quad  n\geq n_{_{0}}.$$
Consequently, for all $n,m\geq n_{_{0}}$, it follows that
$$\eqalign {\Vert x_{_{n}}-x_{_{m}}\Vert ^{\alpha }_{_{E(\nm)}}
&\leq \Vert e(x_{_{n}}-x_{_{m}})e\Vert ^{\alpha }_{_{E(\nm)}}
+3^{1/s^{\prime }}\left (\Vert x_{_{n}}\Vert ^k_{_{E(\nm)}}
-\Vert ex_{_{n}}e\Vert ^k_{_{E(\nm)}}\right )^{1/s}\cr
&+3^{1/s^ {\prime }}\left (\Vert x_{_{m}}\Vert ^k_{_{E(\nm)}}
-\Vert ex_{_{m}}e\Vert ^k_{_{E(\nm)}}\right )^{1/s} \cr
&\leq \Vert e(x_{_{n}}-x_{_{m}})e\Vert ^{\alpha }_{_{E(\nm)}}
+2.3^{1/s^{\prime }}(1-(1-\delta )^k)^{1/s}.\cr }$$
We obtain that
$$\inf _{_{n\neq m}}\Vert x_{_{n}}-x_{_{m}}\Vert _{_{E(\nm)}}
\leq Q(1-(1-\delta )^k)^{1/k},$$
where $Q=2^{1/\alpha }3^{ 1/(\alpha s^{\prime }) }$. Thus, given $\epsilon
\in
(0,Q)$, we may take $\delta (\epsilon) $ 
to be
given by the formula
$$ \delta (\epsilon )=1-\left (1-(\epsilon /Q)^k\right )^{1/k} \ ;$$
 and by this the proof of the theorem is complete.

In order to state some explicit consequences of the preceding theorem,
let us
recall (cf. [LT] p.142 ) that, for $0<p\leq q<\infty $, the Lorentz
space
$L_{_{q,p}}(\r)$ consists of all Lebesgue measurable functions $f$
on
$\r$ for which
$$\Vert f\Vert _{_{q,p}}=\left (\int _0^{\infty }\mu
_{_{t}}^p(f)dt^{p/q}\right
)^{1/p}<\infty $$
It is shown in [Cr] Proposition 3.2 that $L_{_{q,p}}(\r)$ is
$p-$convex with $M^{(p)}\left (L_{_{q,p}}(\r)\right )=1 $ and
satisfies a
lower $q-$estimate with constant $1$, in the case that $p\geq 1$;
however, it is
not difficult to see that this restriction is unnecessary and further that
  $L_{_{q,p}}(\r)$ is a $p$-Banach space when $0 < p < 1$.
 We may now state the following consequence of Theorem 2.7.

\noindent {\bf Corollary 2.8.}\quad {\sl If  $0<p\leq q<\infty $ then
the
 Lorentz space $L_{_{q,p}}(\nm)$ has the uniform Kadec-Klee-Huff
property
with respect to local convergence in measure.}

The corollary appears to be new even in the case that $\nm $ is
commutative,
and in this setting, the corollary extends earlier results in [Le2,3].
As
has been noted earlier,  if \ $\nm$ \ is the space \ ${\cal L}(H)$ \ of
all
bounded linear operators on \ $H$ \ and \ $\tau$ \ is the standard
trace, then
convergence locally in measure is precisely convergence for the weak
operator
topology. In this case, the corollary implies that the Schatten
$p-$classes
${\cal C}_{_{p}},0<p<\infty $, have the uniform Kadec-Klee-Huff property
for
the weak operator topology.  For the Schatten class ${\cal C}_{_{1}}$,
this
is again due to Lennard [Le1].  Of course, if $1<p<\infty $, then the
result
for ${\cal C}_{_{p}}$ may be derived directly by using  uniform
convexity and Clarkson's inequalities.

We remark further, that if $\nm $ is commutative and the underlying
measure space is $\sigma-$finite , or if  $\nm$ is the space  ${\cal
L}(H)$
of all bounded linear operators on a separable Hilbert space $H$ with
the usual
trace, or if $\tau (1)<\infty $, then the topology of local convergence
in
measure is a metric topology.  Accordingly, when combined with
[Le2] Corollary 3.3 and Theorem 4.2, our Theorem 2.7 immediately yields
results
concerning normal structure and fixed point properties in each of the
above
cases. We leave precise formulations to the
interested reader.  We state, however a further consequence of
Theorem 2.7 and [BDDL] Theorem 3.1, which seems to be of independent interest.

\noindent {\bf Corollary 2.9}\quad {\sl If $\left (\Omega,\Sigma,\mu
\right
)$ is a measure space, if $1\leq p<\infty $, and if $E$
satisfies the conditions of Theorem 2.7 with $\alpha =1$, then the
Lebesgue-Bochner space $L^p\left (\mu ,E(\nm)\right )$ has the following
property.  For each $\epsilon >0$, there exists $\delta \in (0,1)$ such
that
whenever $\{f_{_{n}}\}$ is a sequence in ${\bf B}_{_{L^p\left
(E(\nm)\right )}}$
and $f\in L^p\left (E(\nm)\right )$ with}
$f_{_{n}}(\omega )\to _{_{n}}f(\omega ) ${\sl locally\ in\ measure\
for\ almost\ all} $\omega \in \Omega \ ${\sl and}
$\Vert f\Vert _{_{L^p\left (E(\nm)\right )}}>1-\delta,$
 {\sl then} {\sl  for some positive integer} $N$
$$\sup_{n,m\geq N}\Vert f_{_{n}}-f_{_{m}}
\Vert _{_{L^p\left (E(\nm)\right)}}\leq \epsilon \ .$$

 We remark that, as noted in [BDDL], the above result is still true if
 $0 < p < 1$. Moreover, the result still holds if  $0 < \alpha < 1$.

\noindent {\bf Theorem 2.10 }\quad {\sl Let} $E$ {\sl be a symmetric Banach
function} {\sl
 space on} $\r$. {\sl If} $E$ {\sl has} $UKKH$ $(lcm)$, {\sl then} $E$
 {\sl satisfies a lower $q$-estimate for some finite} $q$.

\noindent {\bf Proof}\quad If $E$ fails to satisfy a lower $q$-estimate
 for some finite
$q$, then it follows from [LT] 1.f.12 that for every $\epsilon >0$ and
$n\in
\jn$, there exists a (finite) sequence $\{x_{_{i}}\}_{_{i=1}}^n$ of
mutually disjoint elements of $E$ such that
$$\max_{1\leq i\leq n} |a_{_{i}}|\leq \Vert \sum _{i=1}^n
a_{_{i}}x_{_{i}}\Vert
_{_{E}} \leq (1+\epsilon )\max_{1\leq i\leq n} |a_{_{i}}|, $$
for all choices of scalars $ \{a_{_{i}}\}_{_{i=1}}^n.$

Let $\delta>0$ be fixed and choose disjoint elements $x,y \in E$ such
that
$$\max\{|a|,|b|\}\leq \Vert ax+by\Vert _{_{E}}\leq (1+\delta
)\max\{|a|,|b|\}$$
for every choice of scalars $a,b$.  Observe that
$$1\leq \Vert x\Vert _{_{E}},\Vert y\Vert _{_{E}},\Vert x+y\Vert _{_{E}}
\leq 1+\delta .$$
Set $$ x^{\prime }=\Vert x+y\Vert _{_{E}}^{-1}x,\quad
y^{\prime }=\Vert x+y\Vert _{_{E}}^{-1}y$$
so that $\Vert x^{\prime }+y^{\prime }\Vert _{_{E}}
=1 $ and
$$(1+\delta )^{-1}\leq   \Vert x^{\prime }\Vert _{_{E}},
\Vert y^{\prime }\Vert _{_{E}} \le 1+\delta \ .$$
For each $n=0,1,2,\dots ,$ set $I_{_{n}}=[2n+1,2n+2)$ and let $T=\cup
_{n=0}^\infty I_{_{n}}$.  Let $\{T_{_{n}}\}_{_{n=1}}^{\infty }$ be a
 measurable partition of
$T$ into a mutually disjoint sequence of subsets of $T$ with
$|T_{_{n}}|=\infty,
\ n=1,2, \dots .$  Let $\phi :\r \to \r \backslash T$ and $\phi _{_{n}}
:\r \to T_{_{n}},\ n=1,2, \dots $ be measure preserving bijections. Set
$$x^{\prime \prime }=x^{\prime }\circ \phi ^{-1},\ y_{_{n}}=
y^{\prime }\circ \phi ^{-1}_{_{n}},\ z_{_{n}}=x^{\prime \prime
}+y_{_{n}}$$
for $ n=1,2,\dots .$  It follows that
$$\Vert z_{_{n}}\Vert _{_{E}}=\Vert x^{\prime \prime }+y_{_{n}}\Vert
_{_{E}}
=\Vert x^{\prime }+y^{\prime }\Vert _{_{E}} =1,$$
and
$$\eqalign { \Vert z_{_{n}}-z_{_{m}}\Vert _{_{E}}=
\Vert x^{\prime \prime }+y_{_{n}} -( x^{\prime \prime }+y_{_{m}}) \Vert
_{_{E}}&= \Vert y_{_{n}} -y_{_{m}} \Vert
_{_{E}}\cr
&\geq \Vert y_{_{n}}\Vert _{_{E}}=\Vert y^{\prime }\Vert _{_{E}}\geq
(1+\delta )^{-1}\cr }$$
for all $n,m= 1,2, \dots $ with $ n \ne m $.
  Since $y_{_{n}}\to 0$ locally in measure,
it
follows that $z_{_{n}}\to x^{\prime \prime } $ locally in measure.
However,
since $\Vert x ^{\prime \prime } \Vert _{_{E}}\geq (1+\delta )^{-1} $
and since
$\delta >0$ is arbitrary, it follows that $E$ does not have $UKKH(lcm)$.

If $E$ is a symmetric Banach function space on $\r$, and if $E$
satisfies a lower
q-estimate for some $1 <q<\infty $, then it follows from [LT] 1.f.7 that
$E$ is
$r-$concave for every $r>q$. For each such $r>q$, the space $E$ may be
renormed
equivalently so that $E$, endowed with the new norm and the same order
is a
symmetric Banach function space on $\r$ which is $r-$concave with
$r$-concavity constant $1$.  See, for example, [LT] 1.d.8, 
1.f.11. We obtain
therefore the following consequences.

\noindent {\bf Corollary 2.11 }\quad {\sl If $E$ is a symmetric Banach
function space on $\r$, then $E$ satisfies a lower $q-$estimate for some
$1<q<\infty $ if and only if there is an equivalent symmetric norm $\Vert
\cdot
\Vert _{_{0}}$ on $E$ such that $(E,\Vert \cdot \Vert _{_{0}})$ has
$UKKH(lcm)$.}

\noindent {\bf Corollary 2.12 }\quad {\sl If $E$ is a symmetric Banach
function space on $\r$ and if $E$ has $UKKH(lcm)$, then there is an
equivalent
norm $\Vert \cdot \Vert _{_{0}}$ on $E(\nm )$ such that $\left (E(\nm ),
\Vert \cdot \Vert _{_{0}}\right )$ has} $UKKH(lcm)$.

\noindent {\bf 3. Final remarks and related results}\quad There are
several variants of the commutative specialization of Theorem 2.7. In
this
section, we formulate one such variant in the setting of Banach
lattices.  We
begin with the following simple characterization of Banach lattices
which are
$q-$concave, for some finite $q$.

\noindent {\bf Proposition 3.1}\quad {\sl If $E$ is a Banach lattice,
then the
following statements are equivalent.\hfil \break
\noindent \quad {\rm (i)} $E$ satisfies a lower $q$-estimate for some
$q<\infty
$.\hfil \break
\noindent \quad {\rm (ii)} $E$ admits an equivalent norm $\Vert \cdot
\Vert
_{_{0}}$ such that $E$ equipped with $\Vert \cdot \Vert_{_{0}}$ and the
same
order is a Banach lattice with the following property: For every $c>0$,
there exists $0<\delta=\delta (c)$ such that, whenever $x,y,z \in E^+$
satisfy
$$x=y+z,\quad y\wedge z=0,\quad \Vert y\Vert _{_{0}}=1,\quad  \Vert
z\Vert _{_{0}}\geq c,$$
it follows that}  $\Vert x\Vert _{_{0}}\geq 1+\delta (c).$

That (i) implies (ii) is a simple consequence of [LT] 1.f.7, 1.d.8.
The
reverse implication follows from [LT] 1.f.12, and the details are
omitted.

 A Banach lattice $E$ which
satisfies the
stated condition (ii) in Proposition 3.1 preceding in its given norm is
said to
have property $(C)$. We remark that it follows from Proposition 3.1 that
if
a Banach lattice $E$ has property $(C)$, then $E$ necessarily has order
continuous norm.

\noindent {\bf Definition 3.2}\quad {\sl Let $E$ be a Banach lattice.
The
sequence $\{x_{_{n}}\} \subseteq E$ is said to converge locally to $x\in
E$ if
and only if for every band projection $Q$ on $E$ there exists a band
projection
$P$ with $0\neq P\leq Q $ and a subsequence $\{y_{_{n}}\}\subseteq
\{x_{_{n}}\}$
such that }$\Vert P(y_{_{n}}-x)\Vert \to _{_{n}}0.$

Let us mention immediately that if $E$ is any Banach function space with
order
continuous norm on some $\sigma-$finite measure space, and if the
sequence
$\{x_{_{n}}\} \subseteq E$ converges to $x\in E$ in measure on every
subset of
finite measure, then it follows from Egorov's theorem that the sequence
$\{x_{_{n}}\}$ converges locally to $x$.

The theorem which follows, and its proof, is an adaptation of a similar
result
due to van Dulst and de Valk [DV], Proposition 3, in the setting of
Banach
spaces with a Schauder basis. We include the details of proof for the
sake of
completeness.

\noindent {\bf Theorem 3.3}\quad {\sl Let $E$ be a Banach lattice.  If
$E$ has
property $(C)$, then $E$ has the following property. For every $\epsilon
>0$,
there exists $\delta=\delta (\epsilon) \in (0,1)$ such that whenever
$\{x_{_{n}}\}$ is a sequence in the unit ball of $ E$  which converges
locally
to $x\in E $, and which satisfies $\Vert x_{_{n}}-x_{_{m}}\Vert >
\epsilon
 \ ,n\neq m ,$ it follows that }$\Vert x \Vert \leq 1-\delta .$

\noindent {\bf Proof }\quad  Let $\epsilon >0 $ be given and assume that
there is no $\delta >0$ satisfying the assertion of the theorem for this
$\epsilon $. For each $c>0$, let $\delta (c)>0 $ be the largest $\delta
$
which satisfies the assertion of property $(C)$, and let $0<\alpha <1$
satisfy $\alpha (1+\delta (\epsilon/2))>1. $  We assume that there
exists a
sequence $\{x_{_{n}}\}$  in the unit ball of $ E$   which satisfies
$\Vert
x_{_{n}}-x_{_{m}}\Vert > \epsilon \ , \ n\neq m $, and which converges
locally
to some $x\in E $ with $\Vert x\Vert >\alpha .$  Observe that local
convergence of the sequence $\{x_{_{n}}\}$ to $x$ together with
order continuity of the norm on $E$ implies that there exists a band
projection $P$ and a subsequence $\{y_{_{n}}\}\subseteq \{x_{_{n}}\}$
such that $\Vert Px\Vert >\alpha $ and $\Vert P(y_{_{n}}-x) \Vert \to
_n 0.$  Consequently there exists a natural number $n_{_{0}}$ such that
$$ \Vert Py_{_{n}}\Vert >\alpha ,\quad \Vert
(I-P)(y_{_{n}}-y_{_{m}}\Vert > \epsilon $$
for all natural numbers $n,m \geq n_{_{0}}, n\neq m.$  For at least one
of the choices $n=n_{_{0}}, n_{_{0}}+1$, it follows that
$$ \alpha < \Vert Py_{_{n}}\Vert \leq 1,\quad \Vert
(I-P)y_{_{n}}\Vert >\epsilon /2.$$ With this choice of index,
set $$w={ y_{_{n}}\over \Vert Py_{_{n}}\Vert },\quad y=P\left
({y_{_{n}}\over \Vert Py_{_{n}}\Vert }\right ),\quad z=(I-P)\left
({y_{_{n}}\over \Vert Py_{_{n}}\Vert }\right ),$$
and observe that $$|w|=|y|+|z|,\quad |y|\wedge |z|=0,\quad \Vert
\ |y|\ \Vert =1,\quad \Vert \ |z|\ \Vert \geq \epsilon /2.$$
Property $(C)$ now implies that $\Vert |w|\Vert \geq 1+\delta
(\epsilon/2))$ so
that
$$\Vert y_{_{n}}\Vert =\Vert |y_{_{n}}|\Vert >\alpha (1+\delta (\epsilon
/2)>1,$$
and this is a contradiction.

\bigskip

\noindent
 {\bf  References }\hfil\break
\smallskip
{\noindent
[Ar]  J. Arazy,  {\it More on convergence in unitary matrix
spaces}, Proc. Amer. Math. Soc.{\bf 83} (1981),44-48.\hfil\break
[BDDL]  M. Besbes, S.J. Dilworth, P.N. Dowling, and C.J. Lennard,{\it New
Convexity and
fixed point properties in Hardy and Lebesgue-Bochner spaces},
J. Funct. Anal.
(to appear)\hfil \break
 [CDLT]  N.L. Carothers, S.J. Dilworth, C.J. Lennard and D.A.
Trautman, {\it A fixed point property for the Lorentz space
$L_{p,1}(\mu)$},
Indiana Univ. Math. J.{\bf 40} (1991) 345-352.\hfil \break
[Cr]  J. Creekmore, {\it Type and cotype in Lorentz $L_{_{pq}}$ spaces},
Indag. Math., {\bf 43} (1981)145-152.\hfil\break
[CS]  V.I. Chi'lin and F.A. Sukochev, {\it Convergence in measure in
admissible non-commut-ative symmetric spaces}, Izv. Vys\v s U\v
ceb.Zaved {\bf 9}(1990) 63-90(Russian).\hfil\break
[CKS]  V.I. Chi'lin, A.V. Krygin
and P.A. Sukochev, {\it Uniform and local uniform convexity of symmetric
spaces of measurable operators,}, Dep VINITI N5620-B90(1990) 24pp.
(Russian).\hfil\break
[DH] S.J. Dilworth and Y.P. Hsu, {\it The uniform Kadec-Klee property
in the trace ideals ${\cal C}_E$,} preprint 1993. \hfil \break
[DDP1] P.G. Dodds, T.K. Dodds and B. de
Pagter, {\it Non-commutative Banach function
  spaces}, Math. Z. {\bf 201} (1989), 583-597.\hfil\break
[DDP2]  P.G. Dodds, T.K. Dodds and B. de Pagter, {\it
Non-commutative  K\"othe duality},Trans. Amer. Math. Soc. (to
appear).\hfil\break
 [DL1] P.N. Dowling and C.J. Lennard, {\it Kadec-Klee
properties
of vector-valued Hardy spaces}, Math. Proc. Camb. Phil. Soc. {\bf 111} (1992)
 535-544.\hfil
\break
 [DL2] P.N. Dowling and C.J. Lennard, {\it On uniformly H-convex complex
quasi-Banach spaces}, Bull. des Sci. Math. (to appear) \hfil \break
[DV] D.van Dulst and V.de Valk,  {\it (KK)-properties, normal structure
and fixed points of nonexpansive
mappings in Orlicz sequence spaces,} Can. J. Math. {\bf 38} (1986)
728-750.\hfil
\break
[Fa]  T. Fack, {\it Sur la notion de valeur
caract\'eristique}, J. Operator Theory {\bf 7} (1982),
307-333.\hfil\break
[FK]  T. Fack and H. Kosaki, {\it Generalized
s-numbers of $\tau$-measurable  operators}, Pacific J.
  Math. {\bf 123} (1986), 269-300.\hfil\break
[Fr]D.H. Fremlin,{\it Stable subspaces of $L^1+L^{\infty }$}, Math. Proc. Camb.
Philos. Soc. {\bf 64} (1968) 625-643.\hfil \break
 [GK]  I.C. Gohberg and M.G. Krein, {\it Introduction to the theory of
non-selfadjoint operators},
  Translations of Mathematical Monographs, vol.18, AMS
(1969).\hfil\break
[H] R. Huff, {\it Banach spaces which are nearly uniformly convex}, Rocky
 Mountain J. Math. {\bf 10} (1980) 743-749. \hfil \break
[KPS]  S.G. Krein, Ju.I. Petunin and E.M. Semenov, {\it Interpolation of
linear operators},
  Translations of Mathematical Monographs, vol.54, AMS
(1982).\hfil\break
[Le1]  C. Lennard, {\it $ {\cal C}_1$ is uniformly Kadec-Klee}, Proc.
Amer. Math. Soc. {\bf 109} (1990) 71-77.\hfil \break
[Le2] C. Lennard, {\it A new convexity property that implies a fixed
point
property for $L_1$, } Studia Math. {\bf 100} (1991) 95-108.\hfil \break
[Le3] C. Lennard, {\it Operators and geometry of Banach spaces}, PhD
Thesis, Kent State University, 1988.\hfil \break
[LT]  J. Lindenstrauss and L. Tzafriri, {\it Classical}
 {\it Banach} {\it Spaces II},
 {\it Function} {\it Spaces}, Spring-er-Verlag, 1979.\hfil \break
[M] C.A. McCarthy, {\it $c_p$}, Israel J. Math. {\bf 5} (1967), 249-271.
\hfil\break
[Ne]  E. Nelson, {\it Notes on non-commutative integration}, J.
Funct. Anal. {\bf15} (1974), 103-116.\hfil \break
[Ov]  V.I. Ovcinnikov, {\it s-numbers of measurable operators},
Funktsional'nyi
Analiz i Ego Prilozheniya {\bf 4} (1970) 78-85 (Russian).\hfil \break
[Si]  B. Simon, {\it Convergence in trace ideals}, Proc. Amer. Math. Soc.
{\bf 83}(1981), 39-43. \hfil \break
[Su] F.A. Sukochev, {\it On the uniform Kadec-Klee property,} International
Conference dedicated to the 100-th birthday of S.Banach, Lwow, 1992.\hfil
\break
 [Te]  M. Terp, {\it $L^p$-spaces associated with von Neumann algebras,}
Notes, Copenhagen University (1981).\hfil \break 
[X]  Q. Xu, {\it Analytic
functions with values in lattices and symmetric spaces of measurable operators,}
Math. Proc. Camb. Phil. Soc. {\bf 109 } (1991) 541-563.\hfil \break
 }

\noindent
 {\bf Addresses }\hfil\break
\smallskip

\settabs 2 \columns
         \+P.G. Dodds&T.K. Dodds\cr
         \+School of Information Sci. $\&$ Technology&School of
  Information Sci. $\&$ Technology\cr
         \+The Flinders University of South Australia&The
  Flinders University of South Australia\cr
         \+GPO Box 2100, Adelaide 5001&GPO Box 2100, Adelaide 5001\cr
         \+South Australia&South Australia\cr
         \bigskip
         \+P.N. Dowling&C.J. Lennard\cr
         \+Department of Mathematics and Statistics&Department of Mathematics
and Statistics\cr
         \+Miami University&University of Pittsburgh\cr
         \+Oxford, Ohio 45056&Pittsburgh, Pennsylvania 15260\cr
         \bigskip
         \+F.A. Sukochev\cr
         \+Department of Mathematics\cr
         \+Tashkent State University\cr
         \+Vuzgorodok, 700095, Tashkent\cr
         \+Uzbekistan\cr

\end